\documentclass[11pt]{article}
 \usepackage{benstyle}
 
\title{A note on compressed sensing of structured sparse wavelet coefficients from subsampled Fourier measurements}
\author{Ben Adcock\footnote{Department of Mathematics, Purdue University, 150 N. University Street, West Lafayette, IN 47907, USA} \and Anders C. Hansen \footnote{DAMTP, Centre for Mathematical Sciences, University of Cambridge, Cambridge CB3 0WA, UK} \and Bogdan Roman \footnote{DAMTP, Centre for Mathematical Sciences, University of Cambridge, Cambridge CB3 0WA, UK} }
\date{}

\begin{document}

\maketitle

\begin{abstract}
This note complements the paper \textit{The quest for optimal sampling: 
Computationally efficient, structure-exploiting measurements for compressed 
sensing} \cite{AHRBerlinBookChpt}. Its purpose is to present a proof of a 
result stated 
therein concerning the recovery via compressed sensing of a signal that has 
structured sparsity in a Haar wavelet basis when sampled using a 
multilevel-subsampled discrete Fourier transform.  In doing so, it provides a 
simple exposition of the proof in the case of Haar wavelets and discrete 
Fourier samples of more general result recently provided in \textit{Breaking 
the coeherence barrier: A new theory for compressed sensing} 
\cite{AHPRBreaking}.
\end{abstract}

\section{Introduction}

In many applications of compressed sensing, the image or signal $x \in \bbC^n$ to be recovered is sparse or compressible in an orthonormal wavelet basis $\Phi \in \bbC^{n \times}$.  However, it is well known that the coefficients $c = \Phi^* x$ in such a basis possess far more than mere sparsity.  In fact, they are highly structured: if the vector $c$ of wavelet coefficients is divided into dyadic scales, there is far more sparsity at the finer scales than at the coarser scales.  In \cite{AHRBerlinBookChpt} it was argued that, in order to obtain a better reconstruction with compressed sensing, one should exploit such structure by taking appropriate measurements.  This can be achieved by subsampling the discrete Fourier transform in an appropriate way.  Not only does this lead to improved reconstructions over standard (sub)Gaussian random measurements, it also explains the success of compressed sensing in applications where the measurements naturally arise from the Fourier transform, e.g.\ MRI, X-ray CT, etc.

In this note we provide a short, expositional proof of the corresponding recovery result stated in \cite{AHRBerlinBookChpt} for the case of one-dimensional discrete Fourier measurements with Haar wavelets.  We refer to \cite{AHPRBreaking} for the proof of the corresponding result for general wavelets in the infinite-dimensional setting.  Throughout, we use the same notation as in \cite{AHRBerlinBookChpt}.

\section{Preliminaries}
Let $x = \{ x(t) \}^{n-1}_{t=0} \in \bbC^n$ be a signal.  Denote the Fourier transform of $x$ by
\bes{
\cF x(\omega) = \frac{1}{\sqrt{n}}\sum^{n}_{t=1} x(t) \E^{2 \pi \I \omega t / n},\quad \omega \in \bbR,
}
and write $F \in \bbC^{n \times n}$ for the corresponding matrix, so that
\bes{
F x = \{ \cF x (\omega) \}^{n/2}_{\omega=-n/2+1}.
}
The concern of this note is the recovery of $x$ from a small subset $y \in \bbC^m$ of the measurements $F x$.  We do so using techniques of compressed sensing, by assuming that $x$ is compressible in a Haar wavelet basis.  Let $n = 2^r$ for some $r \in \bbN$.  The Haar basis consists of the functions $\{ \psi \} \cup \{ \phi_{j,p} : j=0,\ldots,r-1,\ p=0,\ldots,2^j - 1 \}$ where
\bes{
\psi(t) = 2^{-r/2},\quad 0 \leq t < 2^r,
}
and 
\bes{
\phi_{j,p}(t) = \left \{ \begin{array}{cl} 2^{\frac{j-r}{2}} & p 2^{r-j} \leq t < (p+\frac12) 2^{r-j}  \\ -2^{\frac{j-r}{2}} & (p+\frac12) 2^{r-j} \leq t < (p+1) 2^{r-j} \\ 0 & \mbox{otherwise} \end{array} \right . .
}
Write $\Phi \in \bbC^{n \times n}$ for the matrix corresponding to this basis, and let $c = \Phi^* x \in \bbC^n$ be the vector of coefficients of $x$.  We divide $c$ into $r$ levels corresponding to wavelet scales:
\bes{
c = (c^{(0)} | \ldots | c^{(r-1)} )^{\top},
}
(note that we now index over $0,\ldots,r-1$, as opposed to $1,\ldots,r$ as was done in \cite{AHRBerlinBookChpt}) where
\bes{
c^{(0)} = \left ( \ip{x}{\psi} , \ip{x}{\phi_{0,0}} \right )^\top \in \bbC^2,
}
and
\bes{
c^{(j)} = \left ( \ip{x}{\phi_{j,0}},\ldots,\ip{x}{\phi_{j,2^{j-1}}} \right )^\top \in \bbC^{2^j}.
}
Let $M_0 = 0$ and
\bes{
M_j = 2^j,\quad j=1,\ldots,r,
}
so that $c^{(j)}$ corresponds to the segment of the vector $c$ with indices $\{ M_{j} + 1,\ldots,M_{j+1} \}$.

We now wish to specify how to subsample the Fourier transform $F x$.  Recall that $F x$ is indexed over $\{ -n/2+1,\ldots,n/2\}$.  Proceeding as in \cite{AHRBerlinBookChpt}, we divide this set up into $r$ frequency bands.  Let
\bes{
W_0 = \{ 0,1\},
}
and
\be{
\label{W_def}
W_{j} = \{ - 2^{j}+1,\ldots,-2^{j-1} \} \cup \{ 2^{j-1}+1,\ldots,2^{j} \},\quad j=1,\ldots,r-1,
}
and note that $W_0,\ldots,W_{r-1}$ form a disjoint partition of $\{ -n/2+1,\ldots,n/2\}$.  Observe that
\bes{
|W_0| = 2,\qquad |W_j| = 2^j,\quad j=1,\ldots,r-1.
}
For $j=0,\ldots,r-1$, we now choose the index set $\Omega_j \subseteq W_j$ uniformly at random of size $| \Omega_j | = m_j$.  If
\be{
\label{Omega_def}
\Omega = \Omega_0 \cup \cdots \cup \Omega_{r-1},\qquad | \Omega | = m = m_0 + \ldots + m_{r-1},
}
then the vector of measurements is given by $y = P_{\Omega} F x$, where the matrix $P_{\Omega} \in \bbC^{m \times n}$ picks out the elements of $F x$ with entries in $\Omega$.  Equivalently, the measurement matrix $A = P_{\Omega} F$ (see \cite{AHRBerlinBookChpt}). 

\rem{
Throughout this note, we shall use the notations $a \lesssim b$ and $a \gtrsim b$ to mean that there exists a constant $C$ independent of all relevant parameters such that $a \leq C b$ or $a \geq C b$ respectively.
}

\section{Main theorem}

Our concern is signals $x$ for which the vector $c$ is not just approximately sparse, but has a distinct sparsity structure within its wavelet scale.  Given the parameters $\mathbf{M} = (M_0,\ldots,M_{r-1})$, we recall from \cite{AHRBerlinBookChpt} that $c$ is $(\mathbf{k},\mathbf{M})$-sparse in levels, where $\mathbf{k} = (k_0,\ldots,k_{r-1}) \in \bbN^r$ if
\bes{
\| c^{(j)} \|_{0} \leq k_j,\quad j=0,\ldots,r-1.
}
If $\Sigma_{\mathbf{k},\mathbf{M}}$ denotes the set of such vectors, then we define the best $(\mathbf{k},\mathbf{M})$-term approximation of an arbitrary $c \in \bbC^n$ by
\be{
\label{sigma_s_m}
\sigma_{\mathbf{k},\mathbf{M}}(c)_1 = \min_{z \in \Sigma_{\mathbf{k},\mathbf{M}} } \| c - z \|_1.
}
In order to recover such an $x$ from noisy measurements $y = A x + e$ with $\| e \|_2 \leq \eta$, we consider the convex optimization problem
\be{
\label{CS_general}
\min_{z \in \bbC^n} \| \Phi z \|_1\quad \mbox{s.t.}\quad \| y - A z \|_2 \leq \eta.
}
The result we shall prove is the following:
\thm{
\label{t:Fourier_Haar}
Let $x \in \bbC^n$ and $\Omega$ be as in \R{Omega_def}.  Let $\epsilon \in (0,\E^{-1}]$ and suppose that
\be{
\label{m_j_est}
m_j \gtrsim \left ( k_j + \sum^{r-1}_{\substack{l=0 \\ l \neq j}} 2^{-\frac{|j-l|}{2}} k_l \right )\log(\epsilon^{-1}) \log(n),\quad j=0,\ldots,r-1.
}
Then, with probability exceeding $1-k \epsilon$, where $k=k_0+\ldots+k_{r-1}$, any minimizer $\hat{x}$ of \R{CS_general} satisfies
\bes{
\| x - \hat{x} \|_2 \leq C \left ( \eta \sqrt{D} ( 1+E \sqrt{k} ) + \sigma_{\mathbf{k},\mathbf{M}}(\Phi^* x)_1 \right ),
} 
for some constant $C$, where $\sigma_{\mathbf{k},\mathbf{M}}(f)$ is as in 
\R{sigma_s_m}, $
D= 1+ \frac{\sqrt{\log_2\left(6\epsilon^{-1}\right)}}{\log_2(4En\sqrt{k})}$ 
and $E = \max_{j=0,\ldots,r-1} \{ (N_{j}-N_{j-1})/m_j \}$.  If $m_j = |W_j|$, $j=0,\ldots,r-1$, then this holds with probability $1$.
}

We refer to \cite{AHRBerlinBookChpt} for a detailed discussion on the implications of this result.  However, note that \R{m_j_est} asserts that we require near-optimal number of measurements $m_j$ in the $j^{\rth}$ frequency band to recover the $k_j$ significant wavelet coefficients in the corresponding $j^{\rth}$ wavelet band.

\section{Proof of Theorem \ref{t:Fourier_Haar}}

\subsection{Setup}
Let $U =\bbC^{n \times n}$ be given by $U = F \Phi^*$.  There is a natural division of $U$ into blocks defined by the sampling and sparsity bands.  Let $U_{jl}$ be restriction of $U$ to rows with indices in $W_j$ and columns with indices $\{ M_{l}+1,\ldots,M_{l+1} \}$.  Note that the entries of $U_{j,l}$ are
\bes{
(U_{jl})_{\omega,p} = \cF \phi_{l,p}(\omega),\quad \omega \in W_j,\  p =0,\ldots,2^j-1, j=0,\ldots,r-1,l=1,\ldots,r-1,
}
and
\bes{
(U_{j0})_{\omega,0} = \cF \psi (\omega),\quad (U_{j,0})_{\omega,1} = \cF \phi_{0,0}(\omega),\quad \omega \in W_j,j=0,\ldots,r-1.
}

Recall the coherence $\mu(V)$ of matrix $V \in \bbC^{n \times n}$ is defined by $\mu(V) = \max_{j,l=1,\ldots,n} | V_{j,l} |^2$.   As in \cite[Def.\  4]{AHRBerlinBookChpt}, define the $(j,l)^{\rth}$ local coherence of the matrix $U$ by
\be{
\label{local_coh}
\mu(j,l) = \sqrt{\mu(U_{jl})} \max_{l'=0,\ldots,r-1} \sqrt{\mu(U_{jl'})}
}
Note that the second term is the coherence of the $2^j \times 2^r$ submatrix of $U$ formed by concatenating only those rows in $W_j$.  Given a vector $\mathbf{k} = (k_0,\ldots,k_{r-1})$, we also define the relative sparsities (see \cite[Def.\ 5]{AHRBerlinBookChpt}) by
\be{
\label{Kj_def}
K_j = \max_{\substack{z \in \Sigma_{\mathbf{k},\mathbf{M}} \\ \| z \|_{\infty} \leq 1}} \nm{ \sum^{r-1}_{l=0} U_{jl} z^{(l)} }^2_2.
}
With these definitions in hand, \cite[Thm.\ 1]{AHRBerlinBookChpt} gives that the conclusions of Theorem \ref{t:Fourier_Haar} hold, provided $m_0,\ldots,m_{r-1}$ satisfy the following two conditions:
\begin{itemize}
\item[$(i)$] We have
\be{
\label{ineq_1}
m_j \gtrsim |W_j| \left ( \sum^{r-1}_{l=0} \mu(j,l) k_l \right ) \log (\epsilon^{-1}) \log(n),\quad j=0,\ldots,r-1.
}
\item[$(ii)$] For all $\tilde{k}_0,\ldots,\tilde{k}_{r-1} \in (0,\infty)$ satisfying
\bes{
\tilde{k}_0+\ldots+\tilde{k}_{r-1} \leq k_0 + \ldots + k_{r-1},\quad \tilde{k}_{j} \leq K_j,
}
we have $m_j \gtrsim \tilde{m}_j \log(\epsilon^{-1}) \log(n)$, where $\tilde{m}_j$ satisfies
\be{
\label{ineq_2}
1 \gtrsim \sum^{r-1}_{j=0} \left ( \frac{|W_j|}{\tilde{m}_j} - 1 \right ) \mu(j,l) \tilde{k}_j,\quad l=0,\ldots,r-1.
}
\end{itemize}
Thus, to prove Theorem \ref{t:Fourier_Haar}, we need only show that \R{m_j_est} implies \R{ineq_1} and \R{ineq_2}.  To do this, we need to estimate the local coherences $\mu(j,l)$ and the relative sparsities $K_j$.  These are subjects of the next two subsections.

\subsection{The local coherences $\mu(j,l)$}
We commence with the following lemma:

\lem{
\label{l:FT_Haar}
For $\omega \in \{ -2^{r-1}+1,\ldots,2^{r-1} \}$, we have
\bes{
\cF \psi(\omega) = \left \{ \begin{array}{ll} 1 & \omega = 0 \\ 0 & \mbox{otherwise} \end{array} \right .,
}
and 
\bes{
\cF \phi_{j,p}(\omega) = \left \{ \begin{array}{cl} 0 & \omega = 0 \\  2^{j/2-r} \E^{2 \pi \I \omega p /2^j} \frac{\left (1- \E^{2 \pi \I \omega/2^{j+1}} \right )^2}{1-\E^{2 \pi \I \omega / 2^r}} & \mbox{otherwise} \end{array} \right . .
}
}
\prf{
The first statement is trivial.  For the second, we proceed by direct computation:
\eas{
\cF \phi_{j,p}(\omega) & = \frac{2^{\frac{j-r}{2}} }{\sqrt{n}} \sum_{p 2^{r-j} \leq t < (p+1/2) 2^{r-j}} \E^{2 \pi \I \omega t/n} -  \frac{2^{\frac{j-r}{2}} }{\sqrt{n}} \sum_{ (p+1/2) 2^{r-j} \leq t < (p+1) 2^{r-j}} \E^{2 \pi \I \omega t/n} 
\\
& = \frac{2^{\frac{j-r}{2}} }{\sqrt{n}} \E^{2 \pi \I \omega p 2^{r-j} / n} \sum^{2^{r-j-1}-1}_{s=0} \E^{2 \pi \I \omega s /n} - \frac{2^{\frac{j-r}{2}}}{\sqrt{n}} \E^{2 \pi \I \omega (p+1/2) 2^{r-j} / n} \sum^{2^{r-j-1}-1}_{s=0} \E^{2 \pi \I \omega s/n}
\\
& = 2^{j/2-r} \left ( \E^{2 \pi \I \omega p /2^j} -  \E^{2 \pi \I \omega (p+1/2) /2^j}  \right ) \sum^{2^{r-j-1}-1}_{s=0} \E^{2 \pi \I \omega s/n}
\\
& = 2^{j/2-r} \left ( \E^{2 \pi \I \omega p /2^j} -  \E^{2 \pi \I \omega (p+1/2) /2^j}  \right ) \left ( \frac{\E^{2 \pi \I \omega 2^{r-j-1}/n}-1}{\E^{2 \pi \I \omega / n} - 1} \right )
\\
& = 2^{j/2-r} \E^{2 \pi \I \omega p /2^j} \left ( 1 -  \E^{2 \pi \I \omega /2^{j+1}}  \right ) \left  ( \frac{\E^{2 \pi \I \omega/2^{j+1}}-1}{\E^{2 \pi \I \omega / 2^r}-1} \right ),
}
as required.
}

We now have the following:
\lem{
\label{l:local_coh}
The local coherences $\mu(j,l)$ satisfy
\bes{
\mu(j,l) \lesssim 2^{-j} 2^{-|j-l|/2},\quad j,l=0,\ldots,r-1.
}
}
\prf{
Recalling the definition \R{local_coh}, we see that it suffices to show that
\bes{
\mu(U_{jl}) \lesssim 2^{-j} 2^{-|j-l|},\quad j,l=0,\ldots,r-1.
} 
Let $\omega \in W_j$.  Then 
\be{
\label{omega_range}
2^{j-1} \leq | \omega | \leq 2^j.
}
Recall also that
\bes{
|\sin \pi t| \leq \pi |t|,\quad \forall t \in \bbR,\qquad | \sin \pi t | \geq 2 t,\  | t | \leq 1/2.
}
Thus
\bes{
2^{j-r} \leq \left | \sin (\pi \omega / 2^r ) \right | \leq \pi 2^{j-r},\quad \omega \in W_j.
}
Applying this and Lemma \ref{l:FT_Haar} now gives
\be{
\label{Fphi_bound}
| \cF \phi_{l,p}(\omega) | = 2^{l/2-r+1}  \frac{\left | \sin ( \pi \omega / 2^{l+1} ) \right |^2}{| \sin (\pi \omega / 2^r) |} \lesssim 2^{l/2-j} \left | \sin ( \pi \omega / 2^{l+1} ) \right |^2,\quad \omega \neq 0.
}
Recall also that $\cF \phi_{l,p}(0) = 0$.  Suppose now that $l \geq j$.  Then $| \omega | / 2^l \leq 2^{j-l}$ and therefore we get 
\bes{
| \cF \phi_{l,p}(\omega) | \lesssim 2^{-l/2} 2^{j-l} = 2^{-j/2} 2^{-3|j-l|/2},\quad \forall \omega, l \geq j.
}
Conversely, if $l < j$, then we use the fact that $\left | \sin ( \pi \omega / 2^{l+1} ) \right | \leq 1$ to get
\bes{
| \cF \phi_{l,p}(\omega) | \lesssim 2^{l/2-j} = 2^{-j/2} 2^{-|j-l|/2},\quad \forall \omega, l < j.
}
Hence, we find that
\bes{
| \cF \phi_{l,p}(\omega) | \lesssim 2^{-j/2} 2^{-|j-l|/2},\quad \forall \omega,j,l.
}
Since $U_{jl}$ has entries $\cF \phi_{l,p}(\omega)$ for $l \neq 0$, it now follows immediately that
\bes{
\mu(U_{jl}) \lesssim 2^{-j} 2^{-|j-l|},\quad  j=0,\ldots,r-1,l=1,\ldots,r-1.
}
To complete the proof, we need only consider the case $l=0$.  Recall that when $l=0$, the first column of the matrix $U_{j,l}$ has entries $\cF \psi(\omega)$.  However, by Lemma \ref{l:FT_Haar}, $\cF \psi(\omega) =1$ for $\omega = 0 \in W_0$ and $\cF \psi(\omega) = 0$ for $\omega \neq 0$.  Thus $| \cF \psi(\omega) | \lesssim 2^{-j/2} 2^{-|j-0|/2}$.  The second column has entries $\cF \phi_{0,0}(\omega)$, and thus also satisfies the same bound.  Hence we get the case $l=0$ as well.
}

\subsection{The relative sparsities $K_j$}
From the definition \R{Kj_def}, we have
\bes{
\sqrt{K_j} \leq \max_{\substack{z \in \Sigma_{k,M} \\ \| z \|_{\infty} \leq 1}} \sum^{r-1}_{l=0} \| U_{jl} \|_2 \| z^{(l)} \|_2.
}
Note that $\| z^{(l)} \|_2 \leq \sqrt{\| z^{(l)} \|_0} = \sqrt{k_l}$. Hence
\be{
\label{Kl_ineq}
\sqrt{K_j} \leq \sum^{r-1}_{l=0} \| U_{jl} \|_2 \sqrt{k_l},
}
and therefore it suffices to estimate $\| U_{jl} \|_2$.

\lem{
The matrices $U_{jl}$ satisfy
\bes{
\| U_{jl} \|_2 \lesssim 2^{-|j-l|/2},\quad j,l=0,\ldots,r-1.
}
}
\prf{
Suppose that $l=0$ and let $z \in \bbC^2$, $\| z \|_2 = 1$.  Then
\bes{
\| U_{j0} z \|^2_2 = \sum_{\omega \in W_j} \left | \cF \psi(\omega) z_0 + \cF \phi_{0,0}(\omega) z_1 \right |^2 \leq \sum_{\omega \in W_j} \left ( \left | \cF \psi(\omega) \right |^2 + \left | \cF \phi_{0,0}(\omega)  \right |^2 \right ).
}
Recall that $\cF \psi(\omega) = 0$ for $\omega \neq 0$ and $\cF \psi(0) = 1$.  Also $\cF \phi_{0,0}(0) = 0$ and by \R{Fphi_bound} we have $| \cF \phi_{0,0}(\omega) | \leq 2^{-j}$.  Since $| W_0 | = 2 $ and $| W_j | = 2^j$ otherwise, we get $\| U_{j0} z \|^2_2 \lesssim 2^{-j}$.  The result for $l=0$ now follows immediately.  

Suppose now that $l = 1,\ldots,r-1$.  Let $z \in \bbC^{2^l}$, $\| z \|_2 = 1$, and write $g = \sum^{2^l-1}_{p=0} z_p \phi_{l,p}$.  Then
\be{
\label{grizzly_bear}
\| U_{jl}\|^2_2 = \sup_{\substack{z \in \bbC^{2^l} \\ \| z \|_{2}=1}} \sum_{\omega \in W_j} | \cF g(\omega) |^2.
}
By Lemma \ref{l:FT_Haar}, we have $\cF \phi_{l,p}(\omega) = \E^{2 \pi \I \omega p /2^l} \cF \phi_{l,0}(\omega)$.  Hence
\bes{
\cF g(\omega) = \cF \phi_{l,0}(\omega) \sum^{2^l-1}_{p=0} z_p \E^{2 \pi \I \omega p /2^l} =  \cF \phi_{l,0}(\omega) G(\omega/2^l),\qquad G(z) = \sum^{2^l-1}_{p=0} z_p \E^{2 \pi \I p z}.
}
Thus
\ea{
 \sum_{\omega \in W_j} | \cF g(\omega) |^2 &\leq \max_{\omega \in W_j} |\cF \phi_{l,0}(\omega) |^2  \sum_{\omega \in W_j} \left |G(\omega/2^l) \right |^2 \nn
 \\
 &\lesssim 2^{l-2j} | \sin ( \pi \omega / 2^{l+1} ) |^4 \sum_{\omega \in W_j} \left |G(\omega/2^l) \right |^2 \nn
 \\
 & \lesssim \sum_{\omega \in W_j} \left |G(\omega/2^l) \right |^2 \left \{ \begin{array}{ll} 2^{l-2j} & j \geq l \\ 2^{2j-3l} &  j < l  \end{array}\right . , \label{polar_bear}
}
where the second inequality is due to \R{Fphi_bound}.  Since $G(z)$ is periodic with period $1$, we find that 
\be{
\label{black_bear}
\sum_{\omega \in W_j} \left |G(\omega/ 2^l) \right |^2 =  \sum^{2^j-1}_{\omega=0} \left | G(\omega/2^l) \right |^2 .
}
Moreover, since $G$ is a trigonometric polynomial of degree $2^l$, we have
\bes{
\sum^{2^l-1}_{\omega=0} \left | G(\omega/2^l) \right |^2 = 2^l \int^{1}_{0} | G(z) |^2 \D z = 2^{l} \| z \|^2_2 = 2^l.
}
Suppose that $j < l$.  Then by this and \R{black_bear}, we have
\bes{
\sum_{\omega \in W_j} \left |G(\omega/ 2^l) \right |^2 \leq \sum^{2^l-1}_{\omega=0} \left | G(\omega/2^l) \right |^2 = 2^l.
}
Conversely, suppose that $j \geq l$.  By \R{black_bear} and periodicity of $G$, 
\bes{
\sum_{\omega \in W_j} \left |G(\omega/ 2^l) \right |^2 = 2^{j-l} \sum^{2^l-1}_{\omega=0} \left | G(\omega/2^l) \right |^2 = 2^j.
}
Substituting this into \R{polar_bear} and using \R{grizzly_bear} gives
\bes{
\| U_{jl} \|^2_2 \lesssim \left \{ \begin{array}{ll} 2^j 2^{l-2j} & j \geq l \\ 2^l 2^{2j-3l} & j < l \end{array} \right .,
}
and therefore $\| U_{jl} \|^2_2  \lesssim 2^{-|j-l|}$, as required.
}

Using this lemma and \R{Kl_ineq}, we now deduce that
\be{
\label{Kl_ineq2}
K_j \lesssim \left ( \sum^{r-1}_{l=0} 2^{-|j-l|/2} \sqrt{k_l} \right )^2 \lesssim \sum^{r-1}_{l=0} 2^{-|j-l|/2} \sum^{r-1}_{l=0} 2^{-|j-l|/2} k_l \lesssim \sum^{r-1}_{l=0} 2^{-|j-l|/2} k_l .
}

\subsection{Final arguments}
We are now able to complete the proof of the main result, Theorem \ref{t:Fourier_Haar}.  Recall that it suffices to show that \R{m_j_est} implies \R{ineq_1} and \R{ineq_2}.  Consider the right-hand side of \R{ineq_1}.  By Lemma \ref{l:local_coh},
\bes{
|W_j| \left ( \sum^{r-1}_{l=0} \mu(j,l) k_l \right ) \log (\epsilon^{-1}) \log(n) \lesssim\left ( \sum^{r-1}_{l=0} 2^{-|j-l|/2} k_l \right ) \log (\epsilon^{-1}) \log(n).
}
Hence \R{m_j_est} implies \R{ineq_1}.  Similarly, applying Lemma \ref{l:local_coh}  to the right-hand side of \R{ineq_2} gives
\bes{
 \sum^{r-1}_{j=0} \left ( \frac{|W_j|}{\tilde{m}_j} - 1 \right ) \mu(j,l) \tilde{k}_j \lesssim \sum^{r-1}_{j=0}\frac{|W_j|}{\tilde{m}_j} 2^{-j} 2^{-|j-l|/2} \tilde{k}_{j}.
}
Since $|W_j| = 2^j$ and
\bes{
\sum^{r-1}_{j=0} 2^{-|j-l|/2} \lesssim 1,\quad l =0,\ldots,r-1,
}
we see that it suffices to take
\bes{
\tilde{m}_j \gtrsim \tilde{k}_j.
}
By definition, $\tilde{k}_{j} \leq K_j$.  Therefore an application of \R{Kl_ineq2} now gives that \R{m_j_est} implies \R{ineq_2} as well.  This completes the proof of Theorem \ref{t:Fourier_Haar}.

\bibliographystyle{abbrv}
\small
\bibliography{DiscreteFourierHaarRefs}

\end{document}